\documentclass[12pt]{article}
\usepackage{amssymb}
\usepackage{amsfonts}
\usepackage{amsmath}

\newtheorem{theorem}{Theorem}

\newtheorem{lemma}{Lemma}
\newtheorem{proposition}{Proposition}

\newtheorem{definition}{Definition}
\newtheorem{remark}{Remark}
\newenvironment{proof}[1][Proof]{\noindent\textbf{#1.}}

\renewcommand{\theenumi}{(\alph{enumi})}

\begin{document}

\title{Conformal   $\beta$- change in Finsler spaces}%
\author{ S. H. Abed \thanks{Department   of
Mathematics,   Faculty   of   Science,    Cairo University, Cairo,
Egypt.}}
\date{}
\maketitle \noindent{\bf Abstract.}
 We investigate what we call a
conformal $\beta $  - change in Finsler spaces, namely $$
L(x,y)\rightarrow ~^{\ast }L(x,y)=e^{\sigma(x)}L(x,y)+\beta(x,y)$$
where~$\sigma ~$ is a function of $x~ only ~and ~ \beta ( x, y)$ is
a given 1- form.

This change generalizes various types of changes:  conformal
changes, Randers changes and $\beta$ - changes.

 Under this change, we obtain the relationships between
some tensors associated with $(M,L)$ and the corresponding tensors
associated with $(M,{^\ast}L)$.
 We investigate some $\sigma$- invariant tensors . This investigation
allows us to give an answer to the question: Are the properties of
C-reducibility, $S_3$-likeness and $S_4$-likeness invariant under a
conformal $\beta$ - change?

\section{Introduction and Notations}%
 Let $(M,L)$ be a Finsler space, where M is an n-dimensional
differentiable manifold equipped with a fundamental function L.
Given a function $\sigma$, the change $$ L(x,y)\longrightarrow
e^{\sigma (x)}L(x,y)$$ is called a conformal change. The conformal
theory of Finsler spaces has been initiated by M.S. Kneblman
\cite{Knebelman} in 1929 and has been deeply investigated by many
authors: \cite{Hashiguchi}, \cite{Izumi}, \cite{Kitayama},... etc

 In 1941, Randers \cite{Randers} has introduced the Finsler change $$ ~^{r}L(x,y)\longrightarrow  ~^{r}L(x,y)+\beta(x,y)$$ where $^r$L is a
 Riemanian structur and  $\beta$ is a 1-form on M. The resulting space
is a Finsler space. This change has been studied by several authors:
\cite{Matsumoto},\cite{Yasda},...etc.

 The Randers change has been
generalized by Shibata \cite{Shibata} to what is called a $\beta$ -
change
$$ L(x,y)\longrightarrow L(x,y)+\beta(x,y)$$  where L a fundamental Finslerian function. The resulting space
known as a generalized Randers space was studied in \cite{Tamim },
\cite{Kitayama}, \cite{Matsumoto}, \cite{Asuma} and \cite{
Miron},... etc.

 In this paper, we construct a theory which generalizes all the above mentioned changes. In fact, we
consider a change of the form
\begin{equation*}
L(x,y)\longrightarrow \,^{\ast }L(x,y)=e^{\sigma (x)}L(x,y)+\beta (
x, y),
\end{equation*}%
where $\sigma$ is a function of x and $\beta ( x, y ) =
b_{i}(x)y^{i}$ is a 1- form on $M$, which we call  a conformal
$\beta $ - change. This change generalizes various type of changes.
When $\beta = 0$, it reduces to a conformal change. When $\sigma =
0$, it reduces to a $\beta$- change and consequently to a Randers
change.

 We obtain the relationships between some tensors associated with $(M,L)$ ( the fundamental
tensor, the $h(hv)$ - torsion and the third curvature tensor) and
the corresponding tensors associated with $(M,^{\ast }L)\,$.

Under the conformal $\beta $ - change, we investigate  some
$\sigma$- invariant tensors (a tensor K is $\sigma$- invariant if
$~^{\ast}K(x,y)=e^{\sigma}K(x,y)$).

This investigation leads us to find out necessary and /or sufficient
conditions for the properties of C-reducibility,  $S_3$-likeness
and  $S_4$-likeness to be invariant under a conformal $\beta$ -
change  (cf. theorems A, B, and C).

     More investigation and development of this theory will be the object of  forthcoming papers.

     Throughout the present paper, ($x^{i}$) denotes the coordinates of a point of the base manifold M and ($y^{i}$)  the supporting element ($\dot{x}^{i}$).

 We use the following notations:

\thinspace  $l_{i}:= \dot{\partial}_{i}L~=\frac{\partial L}{\partial
y^{i}}:$ the~normalized~supporting~element,

$h_{ij}$ $:= L\dot{\partial}_{i}l_{i}=Ll_{ij}:$ the angular metric
tensor,

$\,\,g_{ij}:= \frac{1}{2}\dot{\partial}_{i}$
$\dot{\partial}_{j}L^{2}:$ the fundamental tensor,

$c_{ijk}:=~ \overset{.}{\partial }_{k\,}(g_{ij}/2)\,:$ the(h) hv
-torsion tensor ,

$c_{i}:= g^{jk}\,c_{ijk}:$ the torsion vector,

$c^{k}:=g^{jk}\,\medskip \allowbreak c_{i}:~,c^{2}=c_{i}c^{i}$,

$S_{hijk}:= c_{ijr}c_{h}\,_{k}^{r}\,-c_{ikr}c_{h}\,_{j}^{r}:$ the
components of the third curvature tensor.

\section{Conformal $\beta $-change}%

 We firstly introduce the following definition

\begin{definition}
A change of Finsler metric defined by
\begin{equation}\label{eq1}
L(x,y)\longrightarrow \,^{\ast }L(x,y)=e^{\sigma (x)}L(x,y)+\beta (
x, y)
\end{equation}
where~$\sigma =\sigma (x)~$is a function of $x ~ and ~\beta
(\medskip x,\medskip y)=b_i(x)y^{i}~$ is a 1- form, will be called
a conformal $\beta$-change.

\end{definition}

This change generalizes various  changes studied by Randers
\cite{Randers}, Matsumuto \cite{Matsumoto},
Shibata.\cite{Shibata}...etc.

We assume that $^{\ast }L(x,y)$ enjoys the same properties possessed by $%
L(x,y)$.

As the Finsler space associated to $L$ is denoted by
$(M,L)$, we denote the Finsler space associated to the conformal $\beta -$%
 change by $(M,^{\ast }L)\,$.

 Throughout the whole paper, the geometric objects associated with $^{\ast}L(x,y)~$will be
asterisked.

\begin{definition}
A geometric object $K$ is said to be $\sigma -$%
invariant if it is invariant, up to a factor $e^{\sigma (x)}~$,
under a conformal $\beta -$change:  $^{\ast }K=e^{\sigma (x)}K .$
\end{definition}
It follows from (1) that
\begin{equation}\label{eq2}
\text{\thinspace }^{\ast }l_{i}(x,y)=e^{\sigma (x)}\,l_{i}(x,y)+
b_{i}(x), \vspace{0.25cm}\ \text{\thinspace
}^{*}l_{ij}(x,y)=e^{\sigma (x)}\,l_{ij}(x,y).
\end{equation}%

The angular metric tensor $h_{ij}$ is given in terms of $h_{ij}$ by
\begin{equation}\label{eq3}
^{\ast }h_{ij}\,= ~^{\ast }L^{\ast }l_{ij\,\,} = ~^{\ast }Le^{\sigma
}l_{ij}~=\tau h_{ij},\vspace{0.15cm}\ \text{\thinspace
}\vspace{0.25cm}\ \tau =e^{\sigma (x)}\frac{^{\ast }L}{L}.\,~\
\end{equation}
Then we have the following

\begin{lemma}
$\frac{h_{ij}}{L}~$is $\sigma -$ invariant under a conformal $\beta
-$change :
\begin{equation}\label{eq4}
\frac{^{\ast }h_{ij}}{\,^{\ast }L}=e^{\sigma \,\,}\frac{h_{ij}}{L}.
\end{equation}
\end{lemma}

As \, $h_{ij}=g_{ij}-l_{i}l_{j}, ~ $ equations
(3) give us a relation between the fundamental tensors $%
g_{ij}~ $and $^{\ast }g_{ij}$:

\begin{equation}\label{eq5}
^{\ast }g_{ij}=\,\,\,\tau (g_{ij}-l_{i}l_{j})+\,_{\,\,}^{\ast
}l_{i}\,_{\,\,}^{\ast }l_{j}.
\end{equation}%

 The relation between the corresponding covariant components is
obtained in the form

\begin{equation}\label{eq6}
^{\ast }g^{ij}=\,\,\,\tau ^{-1}g^{ij}+\mu \,l\,^{i\,}l\,^{j}\vspace{0.15cm}%
-\tau ^{-2}(_{\,\,}l\,^{i}b\,^{j}+l\,^{j}b^{i}),
\end{equation}%
\ where $\mu =(e^{\sigma }Lb^{2}+\beta $ $)/^{\,\,\ast }L\tau ^{2}\,\,,\,\
\,\,\vspace{0.15cm}\vspace{0.15cm}\,\,b^{2}=b_{i}b^{i}\,,\vspace{0.15cm}%
\,\,b^{i}=g^{ij}b_{j}$.

Let us introduce the $\ \pi -$vector field
\begin{equation*}
 \overline{m} =\overline{B}-\frac{\beta }{L^{2}}\overline{\eta \,}\,,\,
\,\;  m^{i}=b^{i}-\frac{\beta }{L}l^{\,i}, \ \   m^{2}=m_{i}\,m^{i}.
\end{equation*}%

\begin{lemma}
For a conformal $\beta -$change which is not conformal (i,e $\beta
\neq 0~),\ \overline{m}\neq 0.$
\end{lemma}

In fact , if $\overline{m}=0$, then $m_{i}\,=0$ for all $i$, and
consequently b$_{i}=\frac{\beta }{L}l_{i}\,$ which implies $\beta
=e^{\psi (x)}L,$   for some function $\psi (x)$.

\begin{lemma}
The(h) hv- torsion tensor $^{\ast }c_{ijk}~$%
associated to $^{\ast }F$ is given by
\begin{equation}\label{eq7}
^{\ast }c_{ijk}=\tau \lbrack c_{ijk}+\frac{1}{2^{\ast }L}h_{ijk}],
\end{equation}%
where
\begin{equation}\label{eq8}
h_{ijk}=h_{ij}m_{k}+h_{jk}m_{i}+h_{ki}m_{j}.
\end{equation}%

From the tensor $^{\ast }c_{ijk},$ we obtain the following important
tensors:
\begin{equation}\label{eq9}
^{\ast }c_{i}\,^{r}\,_{j}=c_{i}\,^{r}\,_{j}+\frac{1}{2^{\ast }L}%
(h_{ij}m^{r}+h_{j}\,^{r}m_{i}+h^{r}\,_{i}m_{j})-\,\tau
^{-1}\,c_{ijs}l^{r}b^{s}-\frac{1}{2^{\ast }L\tau }%
(2m_{i}m_{j}+m^{2}h_{ij})l^{r},
\end{equation}

\begin{equation}\label{eq10}
^{\ast }c_{i}=c_{i}+\frac{n+1}{2^{\ast }L}m_{i}\text,
\end{equation}

\begin{equation}\label{eq11}
^{\ast }c^{k}=\tau ^{-2}[\tau c^{k}-c_{\beta }l^{k}+\frac{n+1}{2^{\ast }L}%
(\tau m^{k}-m^{2}l^{k})],
\end{equation}

\begin{equation}\label{eq12}
^{\ast }c^{2}=\tau ^{-1}[c^{2}+\frac{n+1}{^{\ast }L}\,A_{\beta }],
\end{equation}

\begin{equation}\label{13}
^{\ast }c_{\beta \,}=c_{\beta }+\frac{n+1}{2^{\ast }L}m^{2}
\end{equation}%
where $~A_{\beta }=c_{\beta }+\frac{n+1}{4^{\ast }L}%
m^{2},~~~  c_{\beta }=c_{i}b^{i}\,.$
\end{lemma}
\begin{proof}
\begin{description}
  \item[-] Equation(\ref{eq7}) is deduced from the definition of $c_{ijk}$
together with (\ref{eq5})
  \item [-] Equation(\ref{eq9}) is deduced by rasing the index $ k $ in
(\ref{eq7}), using (\ref{eq6})
  \item [-] Equation(\ref{eq10}) is obtained by contracting the subscript $i$ and
the superscript $r$ in (\ref{eq9})
  \item [-] Equation(\ref{eq11}) follows from (\ref{eq10}) by rasing its
subscript, using (\ref{eq6})
  \item [-] Equation(\ref{eq12}) follows directly from (\ref{eq10}) and
(\ref{eq11}) by contracted multiplication
  \item [-] Equation(\ref{13}) is obtained easily from (\ref{eq10}) and
the definition of $c_{\beta }$ by contracted
multiplication.{$\:\:\:\square$}
\end{description}
\end{proof}
\newpage
\begin{lemma}$\ $

\begin{description}
\item [(a)] The relation between $^{\ast }S_{hijk}~$and\ $S_{hijk}\ $ takes
the form
\begin{equation}\label{eq14}
^{\ast }S_{hijk}=\tau S_{hijk}-\frac{\tau }{2^{\ast }L}\left[h_{ik%
\,}H_{jh}+h_{jh}H_{ik}-h_{hk}H_{ij}-h_{ij}H_{hk}\right],
\end{equation}

where  $\ $%
\begin{equation}\label{eq15}
H_{ij}=c_{i}\,^{r}\,_{j}m_{r}+\frac{1}{2^{\ast }L}m_{i}m_{j}+\frac{1}{%
4^{\ast }L}h_{ij}m^{2}.
\end{equation}

\item [(b)] The $v-$ Ricci tensor $^{\ast }S_{ik}$ is written in the form

\begin{equation}\label{eq16}
^{\ast }S_{ik}=S_{ik}-\frac{1}{2^{\ast }L}[A_{\beta
}h_{ik\,}+(n-3)H_{ik}]
\end{equation}%

\item [(c)] The $v-$ scaler curvature tensor is written in the form

\begin{equation}\label{eq17}
^{\ast }S=\tau ^{-1}[S-\frac{n-2}{^{\ast }L}A_{\beta }]\medskip
\end{equation}
\end{description}
\end{lemma}
\begin{proof}
\begin{description}
  \item [(a)] From equations (\ref{eq7}) and (\ref{eq9}) we have
\begin{eqnarray*}
 ^{\ast }c_{ijr}~^{\ast
}c_{h}\,^{r}\,_{k}&=&\tau \lbrack c_{ijr}+\vspace{0.15cm}\frac{1}{2^{\ast }L}%
h_{ijr}][c_{h}\,^{r}\,_{k}+\frac{1}{2^{\ast }L}(h_{hk}m^{r}+h_{h}%
\,^{r}m_{k}+h^{r}\,_{h}m_{k})\\
&-&\,\tau ^{-1}\,c_{hks}l^{r}b^{s}-\frac{1}{%
2^{\ast }L~\tau }(2m_{h}m_{k}+m^{2}h_{hk})l^{r}]\\
&=&\tau c_{ijr}~c_{h}\,^{r}\,_{k}~~+\frac{\tau }{%
2^{\ast
}L}[(c_{i}\,^{r}\,_{j}~h_{hk}+c_{h}\,^{r}\,_{k}~h_{ij})~m_{r}\\
&+&(c_{ijk~}m_{h}~+~c_{jkh~}m_{i}+~c_{khi~}m_{j}~+~c_{hij~}m_{k})~]\\
&+&\frac{\tau }{4^{\ast }L^{2}}[h_{ij}~h_{hk}~~m^{2}+2h_{hk}~m_{i}~m_{j}+2h_{ij}~m_{h}~m_{k}~+h_{jh}~m_{i}~m_{k}\\
&+&h_{jk} m_{i} m_{h}+h_{ih} m_{j} m_{k}+h_{ik}~m_{j}~m_{h}].
\end{eqnarray*}

Similarly, one can obtain $%
^{\ast }c_{ikr}~^{\ast }c_{h}\,^{r}\,_{j}~~$(by interchange $j~$and
$k$). Hence the result.
\item  [(b)] follows from (a) by contracted multiplication, using (\ref{eq6}).

\item [(c)] is obtained from (b), using (\ref{eq6}) again,  by contracted
multiplication. {$\:\:\:\square$}
\end{description}
\end{proof}

\renewcommand{\theenumi}{\arabic{enumi}}
\begin{remark}

The tensor $H_{ij}~$defined by (\ref{eq15}) has the properties:
\begin{enumerate}
  \item $H_{ij}$ is a symmetric tensor :$\,\,H_{ij}=H_{ji},$
  \item  $H_{ij}$ is an indicatory tensor : $%
H_{ij}\,y^{i}=0=H_{ij}\,y^{j},$
  \item  $g^{ij}H_{ij}$ $=A_{\beta }$.
\end{enumerate}

\end{remark}
\section{Geometrical properties of the conformal $\beta $-change}

\begin{definition}
\cite{MMatsumoto} A Finsler space $(M,L)$ of dimension $n\geq
$~$3~$is\thinspace called a C-reducible space   if \medskip the
h(hv)-torsion tens\medskip or $c_{ijk}\,\,$ has the form
\begin{equation}\label{eq18}
c_{ijk}=h_{ij}M_{k}+h_{kj}M_{i}+h_{ki}M_{j}\,\,,\,\,\,\,\,M_{i}=\frac{c_{i}}{%
n+1}.
\end{equation}%
\end{definition}
Define the tensor
\begin{equation*}
K_{ijk}=[c_{ijk}-(h_{ij}M_{k}+h_{kj} M_{i}+h_{ki}M_{j})]/L.
\end{equation*}%

It is clear that $K_{ijk}~$is a symmetric and indicatory tensor.
Moreover $K_{ijk}$ vanishes if and only if the Finsler space is
C-reducible.

\begin{proposition}
Under a conformal $\beta $-change,\,\, the tensor $K_{ijk}$ is
 $\sigma -$invariant:\
\begin{equation*}
 ^{\ast}K_{ijk}=e^{\sigma }K_{ijk}.
\end{equation*}%
\end{proposition}

\begin{proof}
\, Using Equation (\ref{eq7}) together with the definition of
$K_{ijk},$we get
\begin{eqnarray*}
^{\ast }K_{ijk}&=&[^{\ast }c_{ijk}-(^{\ast }h_{ij}~^{\ast
}M_{k}+^{\ast }h\vspace*{0.1cm}_{ij}~\vspace*{-0.15cm}^{\ast
}M_{i}+^{\ast }h_{ki}~^{\ast }M_{j})]\,/^{\ast }L
\\&=&\tau \lbrack (c_{ijk}+\frac{1}{2^{\ast }L}h_{ijk})-%
\vspace*{0.1cm}(^{\ast }h_{ij}~^{\ast }M_{k}+^{\ast }h_{ij}~\vspace*{-0.15cm}%
^{\ast }M_{i}+^{\ast }h_{ki}~^{\ast }M_{j})]\,/^{\ast }L
\qquad \\&=&\tau \lbrack c_{ijk}+\frac{1}{2^{\ast }L}%
(h_{ij}m_{k}+h_{jk}m_{i}+h_{ki}m_{j}.)-\frac{1}{n~+1}\vspace*{0.1cm}%
(h_{ij}~^{\ast }c_{k}+h_{jk}~\vspace*{-0.15cm}^{\ast
}c_{i}+h_{ki}~^{\ast }c_{j})]\,/^{\ast }L
\\&=&\tau \lbrack c_{ijk}+\frac{1}{n~+1}%
(h_{ij}~c_{k}+h_{jk}c_{i}+h_{ki}c_{j}.)]~/^{\ast }L
\\&=&e^{\sigma }[~c_{ijk}-(h_{ij}M_{k}+\vspace*{0.15cm}h_{kj}%
\vspace*{-0.15cm}M_{i}\vspace*{0.15cm}+h_{ki}M_{j})]\,/L=e^{\sigma
}K_{ijk}~.{\:\:\:\square}
\end{eqnarray*}
\end{proof}

Now, Proposition 1 yields
\begin{theorem}

Under a conformal $\beta $-change $L\longrightarrow ~^{\ast }L,$ the
space$~(M,L)~$is C-reducible if and only if the space $(M,^{\ast
}L)~$is C-reducible.

Consequently the C-reducibility property is invariant under this
change.
\end{theorem}

It shout be noticed that Theorem 4-1 and Corollary 4-1 of Shibata
\cite{Shibata} result from the above Theorem as a very special case.
Some results of Matsumoto \cite{Matsumoto} are also contained in the
above Theorem.

\begin{definition}
 \cite{Ikedo} A Finsler space \,$(M,L)~$ of dimension $ n>4\,$ is called an
$S_{4}-$like space  if the vertical curvature tensor
$S_{hijk}$ has the form%
\begin{equation}\label{eq19}
S_{hijk}=h_{jh}M_{ik}+h_{ik}M_{jh}-\medskip
h_{hk}M_{ij}-h_{ij}M_{hk},
\end{equation}%
where $\ M_{ij}$ is the symmetric and indicatory tensor given by $\
M_{ij}=\frac{1}{n-3}[S_{ij}-\frac{Sh_{ij}}{2(n-2)}].$
\end{definition}
Define the tensor
\begin{equation*}
\eta
_{hijk}=[S_{hijk}-(h_{jh}M_{ik}+h_{ik}M_{jh}-h_{hk}M_{ij}-h_{ij}M_{hk})]/L,
\end{equation*}%
It is clear that $\eta _{hijk}$ vanishes if and only if the manifold $%
(M,L)$ is an $S_{4}-$like manifold.

It is not difficult to prove the following.
\begin{lemma}
The tensor $^{\ast }M_{ij}$~$~$is given ~in ter\vspace*{0.15cm}me%
\vspace*{0.15cm}s of $~M_{ij}~~$by
\begin{equation*}
^{\ast }M_{ij}=M_{ij}~-\frac{1}{2^{\ast }L}H_{ij}.
\end{equation*}
\end{lemma}
 In fact, the result follows from (\ref{eq12}) and (\ref{eq16}).

\begin{proposition}
Under a conformal $\ \beta $-change,\, the tensor\, $\eta _{hijk}$
is $ \sigma-$invariant :
\begin{equation*}
^{\ast }\eta _{hijk}=e^{\sigma }\eta_{hijk}
\end{equation*}
\end{proposition}

\begin{proof}
\, Taking Lemma 4a and Lemma 5 into account, we get
\begin{eqnarray*}
 ^{\ast}L ^{\ast }\eta_{hijk}&=&^{\ast }S_{hijk}-(^{\ast }h_{jh\,}{}^{\ast }M_{ik}+ ^{\ast }h_{ik}{}^{\ast }\,M_{jh}-^{\ast }h_{hk}\,^{\ast
}M_{ij}-^{\ast }h_{ij}{}^{\ast }M_{hk})\,\,\\&=&\tau
S_{hijk}-\frac{\tau }{2^{\ast
}L}[h_{jk\,}H_{ih}+h_{ih}H_{jk}-h_{hk}H_{ij}-h_{ij}H_{hk}]\\
 &-&\tau [ h_{jk\,}(M_{ih}-
\frac{1}{2^{\ast }L}H_{ih})+h_{ik\,}(M_{jh}-\frac{1}{2^{\ast }L}%
H_{jh})\\&-&h_{hk\,}(M_{ij}-\frac{1}{2^{\ast }L}H_{ij})-
h_{ij}~(M_{hk}-\vspace*{0.15cm}\frac{1}{2^{\ast }L}H_{hk})]\\
&=&\tau[S_{hijk}-(h_{jk\,}M_{ih}+h_{ih}\vspace*{0.15cm}%
M_{jk}-h_{hk}M_{ij}-h_{ij}M_{hk})]\\
&=&\tau ~L~\eta _{hijk}=\frac{e^{\sigma \ \,\ast }L}{L}~L~\eta
_{hijk}=e^{\sigma \,\ \ast }L~\eta _{hijk}.
\end{eqnarray*}
Hence the result. {$\square$}\end{proof}

Proposition (2) yields
\begin{theorem}
Under~a~conformal $\beta $-change $L\longrightarrow \,^{\ast }L,$
the \vspace*{0.15cm}space $(M,L)$ is $S_{4}-$like if and only if the
space $(M, ^{\ast }L)$ is an $S_{4}-$like.

Consequently, the $S_{4}-$likeness property is invariant under this
change.
\end{theorem}

The above result generalizes Theorem 4-5 (and its Corollary) of
Shibata \cite{Shibata} .

\begin{definition}
\cite{Ikedo} A Finsler space \,$(M,L)~$ of dimension $ n>3 \,$ is
called an $S_{3}-$like space  if the vertical curvature tensor
$S_{hijk}$ has the form
\begin{equation}\label{eq20}
S_{hijk}=\frac{S}{(n-1)\,\,(n-2)}\,[h_{ik}h_{jh}-h_{ij}h_{hk}]
\end{equation}%
\end{definition}

Define the tensor
\begin{equation*}
\zeta _{hijk}\,\,=[\,\,S_{hijk}-\frac{S}{(n-1)\,\,(n-2)}%
\,(h_{ik}h_{jh}-h_{ij}h_{hk})]\,\,lL
\end{equation*}%
It is clear that ~$\zeta _{hijk~}$ vanishes if and only if the
manifold $(M,L)~$is an $S_{3}-$like manifold

\begin{proposition}
Under a conformal $\beta $-change, the tensor $\zeta _{hijk}$ is $%
~\sigma -$invariant  if and only if\,\, $H_{ij}=\frac{1}{n-1}%
A_{\beta }h_{ij}\,.$
\end{proposition}
\begin{proof}
\,\,Using Equation (\ref{eq7}) together with the definition of $%
K_{ijk}$, we get
\begin{eqnarray*}
^{\ast }L^{\ast }\zeta _{hijk}&=&[^{\ast }S_{hijk}-\frac{^{\ast
}S}{(n-1)\,\,(n-2)}\,(^{\ast }h_{ik}~^{\ast }h_{jh}-^{\ast
}h_{ij}~^{\ast }h_{hk})]\\
&=&[\tau S_{hijk}-\frac{\tau }{2~^{\ast }L}(h_{ik\,}H_{jh}+
h_{jh}H_{ik}-h_{hk}H_{ij}-h_{ij}H_{hk})\\
&-&\frac{\tau }{(n-1)\,\,(n-2)}\,(S-\frac{(n-2)}{^{\ast
}L\,}A_{\beta })(h_{ik}~h_{jh}-h_{ij}~h_{hk})]\\
&=&[\tau (S_{hijk}-\frac{S}{(n-1)\,\,(n-2)}
\,(h_{ik}h_{jh}-h_{ij}h_{hk}))]\\
&-&\frac{\tau }{2~^{\ast }L}[(h _{ik\,}H_{jh}+
h_{jh}H_{ik}-h_{hk}H_{ij}-h_{ij}H_{hk})+\frac{1}{(n-1)\,\,}A_{\beta
}(2h_{ik}~h_{jh}-2h_{ij}~h_{hk})]\\
&=&\tau L~\zeta _{hijk}~-\frac{\tau }{2~^{\ast }L}%
[\,(h_{ik\,}(H_{jh}-\frac{1}{(n-1)} A_{\beta
}h_{jh}))+(h_{jh\,}(H_{ik}-\frac{1}{(n-1)\,\,}A_{\beta }h_{ik}))+\\
&+&(h_{hk\,}(H_{ij}-\frac{1}{(n-1)}A_{\beta
}h_{ij}))+h_{ij}(H_{hk}-\frac{1}{(n-1)\,\,}A_{\beta }h_{hk}))]
\end{eqnarray*}
Now the tensor $\zeta _{hijk}^{\ast }~$is  $\sigma $-invariant\,\,($%
^{\ast }\zeta _{hijk}=e^{\sigma }\zeta _{hijk}$)\,\,if and only
if~all terms of the forms $H_{ij}-\frac{1}{(n-1)}A_{\beta
}h_{ij}$ vanish; that is, if and only if the condition  $H_{ij}=\frac{1}{%
n-1}A_{\beta }h_{ij}\,$ holds. {$\square$}\end{proof} Consequently
we get
\begin{theorem}
Under~a~conformal $\beta $-change $%
L\longrightarrow \,^{\ast }L,$  the following two assertions are
equivalent

\begin{enumerate}
  \item The space $(M,L)$ is $S_{3}-\vspace{0.15cm}$like ,
  \item The spac\vspace{0.15cm}e $(M,^{\ast }L)$ is $S_{3}-$like
\end{enumerate}

if and only if the condition  $H_{ij}=\frac{1}{n-1}A_{\beta }h_{ij}$
holds.

Consequently, the $S_{3}-$ likeness property is invariant under this change if and only if $H_{ij}=\frac{1}{%
n-1}A_{\beta }h_{ij}$

\end{theorem}


\begin{thebibliography}{00}
\bibitem{Hashiguchi}
\author{M. Hashiguchi},  \title{On conformal transformation of \textsc{F}insler metric}, J. Math. Kyoto Univ.. 16(1976) pp. 25-50.
\bibitem{Ikedo}
\author{F. Ikedo},  \title{On \textsc{$S_3$}-and \textsc{$S_4$}-like \textsc{F}insler spaces with the \textsc{T}- tensor of a special form}, Tensor, N.S..35(1981) pp. 345-351.
\bibitem{Izumi}
\author{H. Izumi},  \title{Conformal transformations of \textsc{F}insler spaces \textsc{I} and \textsc{II}}, Tensor, N.S.. 31 and 33(1977 and 1980) pp. 33-41 and 337-359.
\bibitem{Kitayama}
\author{M. Kitayama},  \title{Geometry of transformations of \textsc{F}insler metrics}, Hokkaido University of Education, Kushiro Compus.. Japan(2000)
\bibitem{Knebelman}
\author{M. S. Knebelman},  \title{Conformal geometry of generalized metric spaces}, Proc.nat.Acad. Sci.USA..15(1929) pp. 33-41 and 376-379.
\bibitem{MMatsumoto}
\author{M. Matsumoto},  \title{On \textsc{C} - reducible \textsc{F}insler spaces}, Tensor, N.S.. 24(1972) pp. 29-37.
\bibitem{Matsumoto}
\author{M. Matsumoto},  \title{On \textsc{F}insler spaces with \textsc{R}anders metric and special forms of important tensors}, J. Math. Kyoto Univ.. 14(1974) pp. 477-498.
\bibitem{Miron}
\author{R. Miron},  \title{General  \textsc{R}anders space, Lagrange and \textsc{F}insler geometry}, Ed. by P.L. Antonelli and \textsc{R}.\textsc{M}iron..76(1996) pp.123-140.
\bibitem{Randers}
\author{G. Randers},  \title{On the asymmetrical metric in the four- space of general relativity}, Phys. Rev..2(1941)59 pp. 195-199.
\bibitem{Shibata}
\author{C. Shibata},  \title{On invariant tensors of $\beta$ - change of \textsc{F}insler metrics}, J. Math. Kyoto Univ.. 24(1984) pp. 163-188.
\bibitem{Asuma}
\author{C. Shibata and M. Asuma},  \title{C-conformal invariant and tensors of \textsc{F}insler metrics}, Tensor, N.S.. 52(1993) pp. 76-81.
\bibitem{Yasda}
\author{ C. Shibata and H. Shimada and M. Azuma and H. Yasda},  \title{On \textsc{F}insler spaces with \textsc{R}anders metric}, Tensor, N.S..31(1977) pp. 219-226.
\bibitem{Tamim }
\author{A. A. Tamim and N. L. Youssef},  \title{On generalized Randers manifold}, Algebras, Groups and geometrys..16(1999) pp.115-126.
\end{thebibliography}
\end{document}